\documentclass{gtart_a}
\pdfoutput=1

\title{On the stable equivalence of open books in three-manifolds}

\author{Emmanuel Giroux}
\givenname{Emmanuel}
\surname{Giroux}
\address{\'Ecole Normale Sup\'erieure de Lyon\\69364 Lyon cedex 07\\France}
\email{Emmanuel.Giroux@ens-lyon.fr}
\urladdr{}

\author{Noah Goodman}
\givenname{Noah}
\surname{Goodman}
\address{Massachusetts Institute of Technology\\Cambridge\\Massachusetts 02139\\USA}
\email{ndg@mit.edu}
\urladdr{}

\volumenumber{10}
\issuenumber{}
\publicationyear{2006}
\papernumber{4}
\lognumber{0657}
\startpage{97}
\endpage{114}

\doi{}
\MR{}
\Zbl{}


\received{19 September 2005}
\revised{}
\accepted{28 October 2005}
\published{4 March 2006}
\publishedonline{4 March 2006}
\proposed{Rob Kirby}
\seconded{David Gabai, Peter Ozsv\'ath}
\corresponding{}
\editor{}
\version{}

\arxivreference{math.GT/0509555}
\arxivpassword{}

\begin{asciiabstract}
We show that two open books in a given closed, oriented three-manifold admit
isotopic stabilizations, where the stabilization is made by successive plumbings
with Hopf bands, if and only if their associated plane fields are homologous. 
Since this condition is automatically fulfilled in an integral homology sphere,
the theorem implies a conjecture of J Harer, namely, that any fibered link in 
the three-sphere can be obtained from the unknot by a sequence of plumbings and
deplumbings of Hopf bands. The proof presented here involves contact geometry in
an essential way.
\end{asciiabstract}

\AtBeginDocument{\let\bar\wbar%
\let\tilde\wwtilde\let\wt\wwtilde%
\let\wh\what\let\ol\wwbar}

\makeatletter
\def\cnewtheorem#1[#2]#3{\newtheorem{#1}{#3}
\expandafter\let\csname c@#1\endcsname\c@theorem}

\newtheorem{theorem}{Theorem}
\cnewtheorem{proposition}[theorem]{Proposition}
\cnewtheorem{lemma}[theorem]{Lemma}
\cnewtheorem{corollary}[theorem]{Corollary}
\theoremstyle{definition}
\cnewtheorem{definition}[theorem]{Definition}
\cnewtheorem{example}[theorem]{Example}
\cnewtheorem{remark}[theorem]{Remark}
\newtheorem*{acknowledgments}{Acknowledgments}

\makeatother  

\def\D{\mathbf{D}}
\def\G{\mathbf{G}}
\def\S{\mathbf{S}}
\def\SO{\mathrm{SO}}
\def\ul{\underline}

\let\from\co
\def\res#1{\,\vert\,{}_{#1}}
\def\abs#1{\lvert #1 \rvert}
\def\csum{\mathbin\#}
\def\id{\mathop{\mathrm{id}}}
\def\Int{\mathop{\mathrm{Int}}\nolimits}
\def\PF{\mathcal{PF}}
\def\Hom{\mathop{\mathrm{Hom}}}
\def\class#1{\mathcal{C}_{}^{#1}}
\def\eps{\varepsilon}

\begin{document}

\begin{abstract}
We show that two open books in a given closed, oriented three-manifold admit
isotopic stabilizations, where the stabilization is made by successive plumbings
with Hopf bands, if and only if their associated plane fields are homologous. 
Since this condition is automatically fulfilled in an integral homology sphere,
the theorem implies a conjecture of J~Harer, namely, that any fibered link in 
the three-sphere can be obtained from the unknot by a sequence of plumbings and
deplumbings of Hopf bands. The proof presented here involves contact geometry in
an essential way.
\end{abstract}

\maketitle

Let $M$ be an oriented three-manifold. An \emph{open book} in~$M$ (also called
\emph{open book decomposition} of~$M$) is a pair $(K, \theta)$ consisting of:
\begin{itemize}
\item
a proper one-dimensional submanifold $K$ in $M$;
\item
a fibration $\theta \from M \setminus K \to \S^1 = \R / 2\pi \Z$ which, in some
neighborhood $N = \D^2 \times K$ of $K = \{0\} \times K$, is the normal angular
coordinate.
\end{itemize}
The submanifold $K$ is called the \emph{binding} of the open book while the
closures of the fibers of~$\theta$ are named \emph{pages}. The binding and the 
pages are cooriented by~$\theta$, and hence they are oriented since $M$~is. On 
the other hand, any page~$F$ of an open book $(K,\theta)$ completely determines
$K = \partial F$ and also (though much less evidently, as seen in the work
of Cerf \cite{Ce}, Laudenbach and Blank \cite{LB} and Waldhausen
\cite{Wa}) $\theta$ up to isotopy relative to~$F$.

Around 1920, as a corollary of his results on branched covers and the braiding
of links, Alexander proved the existence of open books in any closed oriented
three-manifold~$M$. On the other hand, given an open book in~$M$, many others
can be constructed by the following plumbing operation. Let $F \subset M$ be a
compact surface with boundary and $C \subset F$ a proper simple arc. We say that
a compact surface $F' \subset M$ is obtained from~$F$ by \emph{$H^\pm$--plumbing
along~$C$} ---~or, more explicitly, by plumbing a positive/negative Hopf band
along~$C$~--- if $F' = F \cup A^\pm$ where $A^\pm$ is an annulus in~$M$ with the
following properties:
\begin{itemize}
\item
the intersection $A^\pm \cap F$ is a tubular neighborhood of $C$ in~$F$;
\item
the core curve of~$A^\pm$ bounds a disk in $M \setminus F$ and the linking
number of the boundary components of $A^\pm$ is equal to~$\pm1$.
\end{itemize}
According to results of J~Stallings~\cite{St} (see \fullref{sec:a}), if $F$ is a page
of an open book $(K,\theta)$ in~$M$ then any surface~$F'$ obtained from~$F$ by
$H^\pm$--plumbing is also a page of an open book $(K',\theta')$ in~$M$. We will
say that the open book $(K',\theta')$ itself is obtained from $(K,\theta)$ by
$H^\pm$--plumbing. A \emph{stabilization} of an open book $(K,\theta)$ is an open
book $(K',\theta')$ that can be obtained from $(K,\theta)$ by finitely many
successive $H^\pm$--plumbings.

If $M$ is closed, any open book $(K,\theta)$ in~$M$ provides a Heegaard spitting
of~$M$: given two antipodal values in~$\S^1$, the two corresponding pages form a
(smooth) closed surface dividing~$M$ into handlebodies. In other words, open
books may be regarded as special Heegaard splittings, namely, those for which
the splitting surface contains a graph whose inclusion in the handlebody on each
side is a homotopy equivalence (a regular neighborhood of this graph is then a
page of the open book, as well as the closure of its complement in the surface).
A well-known theorem by K~Reidemeister and J~Singer shows that any two
Heegaard splittings of a given closed oriented three-manifold admit isotopic
stabilizations, where the stabilization here is made by successive attachings of
trivially embedded one-handles. Furthermore, if an open book $(K',\theta')$ is a
stabilization of another one $(K,\theta)$, the associated Heegard splitting is a
stabilization of the one associated with $(K,\theta)$. It is therefore natural
to ask whether any two open books in a given closed oriented three-manifold have
isotopic stabilizations. To answer this question, we need one more ingredient.

To any open book $(K,\theta)$ in~$M$, we can associate an oriented plane field
$\xi$ on~$M$ in the following way: outside some product neighborhood $N = \D^2
\times K$ of~$K$ in which $\theta$ is the normal angular coordinate, $\xi$~is
just the plane field tangent to the pages, that is, the kernel of~$D\theta$;
inside~$N$, using oriented cylindrical coordinates $(r, \theta, z)$ with $z \in
\S^1 \sqcup \dots \sqcup \S^1$ parameterizing~$K$, we define~$\xi$ by the form
$f(r) \, dz + r^2 d\theta$ where $f \from [0,1] \to [0,1]$ is positive near~$0$
and zero near~$1$. Clearly, the homotopy class of this oriented plane field does
not depend on the choice of $N$ and $r, z, f$: this is an invariant of the open
book that L~Rudolph considered as an ``enhanced Milnor number'' (L~Rudolph was
actually interested in the case $M = \S^3$ where homotopy classes of plane
fields are parameterized by~$\Z$ through the Hopf invariant: see the proof of
\fullref{c:grothendieck}). The main result of this paper is the following:

\begin{theorem} \label{t:equivalence}
Two open books in a closed oriented three-manifold admit isotopic stabilizations
if and only if their associated oriented plane fields are homologous.
\end{theorem}

An oriented hyperplane field on a closed oriented $n$--manifold~$M$ is a section
of the sphere cotangent bundle $ST^*M$. Two hyperplane fields are homologous if
they define equal homology classes in $H_n (ST^*M; \Z)$, or equivalently, if the
curve in~$M$ consisting of points where they coincide with opposite orientations
is nullhomologous (see \fullref{sec:b}).

If $M$ is an integral homology three-sphere, any two plane fields on~$M$ are
homologous. Moreover, an open book in~$M$ is completely determined up to isotopy
by its oriented binding (for this, which again follows from the results of
Cerf \cite{Ce}, Laudenbach and Blank \cite{LB} and Waldhausen \cite{Wa},
it actually suffices that $M$ be a rational homology three-sphere).
Now recall that a link (namely, a closed oriented one-dimensional submanifold) 
in a closed three-manifold~$M$ is a \emph{fibered link} if it is the oriented 
binding of some open book in~$M$. Thus, using the same terminology for fibered 
links as for open books, we get:

\begin{corollary} \label{c:links}
Any two fibered links in an integral homology three-sphere admit isotopic
stabilizations.
\end{corollary}

In the case of the three-sphere itself, the unknot is a fibered link, and so we
get the following result conjectured by Harer~\cite{Ha}:

\begin{corollary} \label{c:harer}
Any fibered link in the three-sphere can be obtained from the unknot by
finitely many plumbings and ``deplumbings'' of Hopf links.
\end{corollary}

These corollaries also admit specific variants for fibered knots, in which the 
$H^\pm$--plumbings leading to the common stabilization can be performed two by
two so as to give a fibered knot at each stage; in other words, the plumbing of
Hopf links is replaced by the plumbing of positive trefoil knots and 
figure-eight knots (see \fullref{sec:a} for definitions):

\begin{corollary} \label{c:knots}
Any two fibered knots in an integral homology three-sphere admit isotopic
stabilizations obtained from each knot by finitely many plumbings of positive 
trefoil knots and figure-eight knots.
\end{corollary}

Following Neumann and Rudolph~\cite{NR}, we can rephrase the above results
concerning the three-sphere in terms of the Grothendieck groups of fibered links
and of fibered knots (see \fullref{sec:a} for the definition):

\begin{corollary} \label{c:grothendieck}
The Grothendieck group of fibered links in the three-sphere is the free Abelian
group of rank~two generated by the positive and the negative Hopf links, $H^+$
and~$H^-$. Similarly, the Grothendieck group of fibered knots in the sphere is 
the free Abelian group of rank~two generated by the positive trefoil knot and 
the figure-eight knot.
\end{corollary}

Let's now say a couple of words about the proof of \fullref{t:equivalence}.
The ``only if'' part is quite easy: an $H^\pm$--plumbing yields an open book that
coincides with the original one in the complement of a ball, and so the homology
class of the associated plane field does not change. To prove the ``if'' part,
we use an invariant of open books more subtle than a homotopy class of plane
fields, namely, an isotopy class of contact structures. The main feature of this
refined invariant is that, according to a previous article by the first
author \cite{Gi}, it determines the open book
up to positive stabilization, that is, stabilization involving only
$H^+$--plumbings (see \fullref{t:pequivalence}). To conclude, we combine this
result with a few observations on the effect of $H^-$--plumbing on our open book
invariants and the classification by Eliashberg~\cite{El} of the so-called
overtwisted contact structures.

\begin{acknowledgments}
The authors thank the American Institute of Mathematics for its support during 
the Fall of 2000. In September 2000, the first author gave a series of lectures 
at Stanford University on the correspondence between contact structures and open
books. The second author was attending these lectures, and both independently 
obtained the results of this paper a few weeks later. 

The first author also thanks the Centre National de la Recherche Scientifique 
for funding his research, and he is very grateful to Alexis Marin for his 
thoughtful comments on Part~B of the text.
\end{acknowledgments}

\section{Plumbing}
\label{sec:a}

Let's first discuss plumbing more carefully, in slightly greater generality (see
the articles by Gabai \cite{Ga1,Ga2}, Harer \cite{Ha} and Stallings
\cite{St}
for further information). For $j \in \{1,2\}$, let $F_j$
be a compact oriented surface in a closed oriented three-manifold~$M_j$ and let
$C_j \subset F_j$ be a proper simple arc. We say that a compact surface $F$ in 
the connected sum $M = M_1 \csum M_2$ is obtained by \emph{plumbing} $F_1$ and
$F_2$ along $C_1$ and~$C_2$ if $F = F_1 \cup F_2$ and $F_1 \cap F_2$ is a square
with median segments $C_1$ and~$C_2$. Thus, the $H^\pm$--plumbing defined in the
introduction is nothing but a plumbing with a positive/negative Hopf band 
in~$\S^3$ (that is, an embedded annulus whose boundary components have linking 
number~$\pm1$) where the arc used in this band connects the two boundary
components.

To see that a surface obtained by plumbing pages of two open books is still a
page of an open book, we start with a simple observation. Consider in~$\R^3$ the
(piece of) open book $(\wh K, \wh\theta)$ whose binding $\wh K$ consists of the
two lines $\{ x=\pm1,\ y=0 \}$ and whose map $\wh\theta \from \R^3 \setminus \wh
K \to \S^1$ is given by
$$ \wh\theta (x,y,z) = \arg \left( \frac {1 + x + iy} {1 - x - iy} \right)
 = \arg ( 1 - x^2 - y^2 + 2iy ) . $$
(Each page of this open book is half of a vertical cylinder containing $\wh K$.)

Let $\wh B$ denote the domain $\{ x^2 + 2y^2 + z^2 \le 2 \}$ and $\wh S$ the
ellipsoid $\partial \wh B$. The map
$$ \wh\rho \from \wh S \longrightarrow \wh S, \quad
   (x,y,z) \longmapsto (z,-y,-x), $$
is an orientation-reversing self-diffeomorphism of order four which permutes the
four points of $\wh S \cap \wh K$ cyclically. Moreover, for $(x,y,z) \in \wh S$,
the identity $y^2 + z^2 - 1 = 1 - x^2 - y^2$ implies that
$$ \wh\theta \circ \wh\rho (x,y,z) = \arg (1-z^2-y^2 - 2iy)
   = \arg (x^2+y^2-1 - 2iy) = \wh\theta (x,y,z) + \pi . $$
Now let $(K_1, \theta_1)$ and $(K_2, \theta_2)$ be open books in closed oriented
three-manifolds $M_1$ and~$M_2$, respectively, and, for $j \in \{1,2\}$, let
$C_j$ be a proper simple arc in the page $\theta_j^{-1}(0) \cup K_j$. Each~$C_j$
has a (big) neighborhood~$W_j$ with an orientation-preserving diffeomorphism
$\phi_j \from W_j \to \R^3$ taking $(K_j \cap W_j, \theta_j \res{ W_j })$ to
$(\wh K,\wh\theta)$ and $C_j$ to the unit segment $\wh C$ of the $x$--axis. Hence
the map
$$ \rho = \phi_2^{-1} \circ \wh\rho \circ \phi_1 \from
   S_1 = \phi_1^{-1}(\wh S) \longrightarrow S_2 = \phi_2^{-1}(\wh S) $$
is an orientation-reversing diffeomorphism and satisfies $\theta_2 \circ \rho =
\theta_1 + \pi$. Therefore, the connected sum
$$ M = M_1 \csum M_2
 = (M_1 \setminus \Int B_1) \mathbin{\mathop\cup\limits_\rho}
   (M_2 \setminus \Int B_2), \qquad
   B_j = \phi_j^{-1} (\wh B), $$
is naturally equipped with an open book $(K,\theta)$ whose binding is the union
$(K_1 \setminus \Int B_1) \cup (K_2 \setminus \Int B_2)$ and whose fibration~$
\theta$ is equal to $\theta_j + (-1)^{j+1} \pi/2$ on $M_j \setminus \Int B_j$, 
$j \in \{1,2\}$. Moreover, the $0$--page of $(K,\theta)$ is (easily seen to be) 
obtained by plumbing the $-\pi/2$--page of $(K_1, \theta_1)$ and the $\pi/2$--page
of $(K_2, \theta_2)$ along the arcs $C_1'$ and~$C_2'$ defined by
\begin{align*}
\phi_1(C_1') & = \bigl\{ x^2 + y^2 = 1,\ y \le 0,\ z = 0 \bigr\}, \\
   \llap{\text{and}\quad}
\phi_2(C_2') & = \bigl\{ x^2 + y^2 = 1,\ y \ge 0,\ z = 0 \bigr\} .
\end{align*}
The open book $(K,\theta)$ is said to be obtained by plumbing $(K_1, \theta_1)$
and $(K_2, \theta_2)$ along $C_1$ and~$C_2$.

In the last section of this paper (see the proof of \fullref{l:homotopy}), we
will call \emph{plumbing ball} for an open book $(K,\theta)$ in~$M$ any ball
$B \subset M$ such that there exists a diffeomorphism $B \to \wh B$ which takes
$(K \cap B, \theta \res B)$ to $(\wh K \cap \wh B, \wh\theta \res{ \wh B })$.

\begin{example} \label{x:plumbing}
Consider $\S^3$ as the unit sphere in $\C^2$:
$$ \S^3 = \left\{ \bigl( r_1 e^{i\theta_1}, r_2 e^{i\theta_2} \bigr) \mid
   r_1^2 + r_2^2 = 1,\ \theta_1, \theta_2 \in \R / 2\pi\Z \right\} . $$
The (unoriented) Hopf link $H = \{ r_1r_2 = 0 \} \subset \S^3$ is the binding of
two open books given by the maps
$$ \theta^\pm \from \S^3 \setminus H \longrightarrow \S^1, \quad
   \bigl( r_1 e^{i\theta_1}, r_2 e^{i\theta_2} \bigr) \longmapsto
   \theta_1 \pm \theta_2 . $$
These two maps orient $H$ in different ways, and we will denote by $H^\pm$ the
Hopf link equipped with the orientation induced by~$\theta^\pm$ (in particular,
the linking number of the components of $H^\pm$ is $\pm1$). On the other hand,
the unknot $U = \{ r_1 = 0 \}$ is the binding of an open book whose fibration is
the map $\theta_1$.

Now let $(K,\theta)$ be an open book in a closed oriented three-manifold~$M$ and
$C$ a proper simple arc in one of its pages. The open book in $M \csum \S^3 = M$
obtained by plumbing $(K,\theta)$ with $(H^\pm, \theta^\pm)$ along~$C$ and an
arc connecting the two components of the Hopf link is what we called earlier the
open book obtained by $H^\pm$--plumbing along~$C$. On the other hand, plumbing
$(K,\theta)$ with~$(U,\theta_1)$ ---~along any arcs~--- always yield an open
book isotopic to $(K,\theta)$.

By plumbing together $(H^+, \theta^+)$ and $(H^\eps, \theta^\eps)$, $\eps \in
\{+,-\}$, along arcs joining the two boundary components, we obtain an open book
in $\S^3$ whose binding is the positive trefoil knot $T^+$ if $\eps = +$ and the 
figure-eight knot~$E$ if $\eps = -$. The plumbing operations with the open books
so obtained will be called \emph{$T^+$--plumbing} and \emph{$E$--plumbing}, 
respectively, provided the arc used in the punctured-torus Seifert surface of 
$T^+$ or~$E$ is non-separating.
\end{example}

The plumbing operation allows us to define a Grothendieck group for fibered 
links in the three-sphere (see the article \cite{NR} by Neumann and
Rudolph). It is the group generated by all (isotopy 
classes of) fibered links in~$\S^3$ in which we impose the relation $[K] = [K']
+ [K'']$ for any triple of fibered links $(K,K',K'')$ such that $K$ is obtained
by plumbing $K'$ and~$K''$ in some way. Clearly, this group is Abelian 
---~because plumbing is a commutative operation~--- and its identity element is
the unknot $U$. The Grothendieck group of fibered knots in the sphere is defined
similarly from the set of all (isotopy classes of) fibered knots in $\S^3$.

\section{Hyperplane fields}
\label{sec:b}

Let $M$ be a connected, oriented $n$--manifold with zero Euler characteristic. We
denote by $\PF(M)$ the (non-empty) space of (co)\,oriented hyperplane fields on
$M$ and by $\PF_c(M)$ its subspace consisting of hyperplane fields that coincide
with a fixed one (arbitrarily chosen) outside of a compact subset of $\Int M$, 
that is, near the boundary and at infinity ---~so $\PF_c(M) = \PF(M)$ if $M$
is closed. This section is a digression in which we investigate the structure of
the set $\pi_0 \PF_c(M)$ of connected components of $\PF_c(M)$ ---~or homotopy 
classes of hyperplane fields. There is nothing new in our discussion, whose key
ideas are due to H~Hopf and L~Pontryagin, but we include it since what we need
is elementary and apparently not so well known (see, however, the articles
by Dufraine \cite{Du}, Gompf \cite{Go} and Turaev \cite{Tu}).

Let $\xi, \eta \in \PF_c(M)$ be hyperplane fields and let $\alpha, \beta$ denote
respective defining one-forms which coincide near the boundary and at infinity.
The first obstruction to the existence of a path joining $\xi$ to~$\eta$ inside
$\PF_c(M)$ is a homology class $c(\xi,\eta)$ in $H_1 (M;\Z) = H_1 ([0,1] \times 
M; \Z)$, namely the class of the zero set of a generic homotopy between $\alpha$
and~$\beta$ with compact support in $[0,1] \times \Int M$. Clearly, for any 
$\xi, \eta, \zeta \in \PF_c(M)$, the following cocycle relations hold:
\begin{align*}
c(\xi,\xi) & = 0 \\
c(\xi,\eta) + c(\eta,\xi) & = 0 \\
c(\xi,\eta) + c(\eta,\zeta) + c(\zeta,\xi) & = 0
\end{align*}
On the other hand, if the linear homotopy $(1-t)\alpha + t\beta$, $t \in [0,1]$,
is generic (that is, transverse to the zero section), then the projection to
$M$ of its zero set is the curve $C(\xi,\eta)$ of points where $\xi$ coincides
with~$-\eta$.

The obstruction class $c(\xi,\eta)$ can also be viewed differently, assuming for
instance that $M$ is closed. As sections of the sphere cotangent bundle $ST^*M$,
the hyperplane fields $\xi$ and~$\eta$ determine homology classes $[\xi], [\eta]
\in H_n (ST^*M;\Z)$. Then consider the long homology exact sequence
$$ \cdots \to H_{n+1} (BT^*M, ST^*M; \Z) \to H_n (ST^*M; \Z) \to H_n (BT^*M; \Z)
   \to \cdots, $$
where $BT^*M$ denotes the ball cotangent bundle. The classes $[\xi], [\eta]$ 
have the same image in $H_n (BT^*M;\Z)$, and so the difference $[\xi]-[\eta]$ is
the image of a class $\tilde c(\xi,\eta) \in H_{n+1} (BT^*M, ST^*M; \Z)$ and 
$c(\xi,\eta) \in H_1(M;\Z)$ is just the intersection of $\tilde c(\xi,\eta)$ 
with the class of the zero section in $H_n (BT^*M;\Z)$. Thus, $c(\xi,\eta) = 0$
if and only if $[\xi] = [\eta]$.

We now need to distinguish homotopy classes of homologous hyperplane fields. Our
tool here is an action of the group $\pi_n \S^{n-1}$ on the set $\pi_0 \PF(M)$.
Take a hyperplane field $\xi \in \PF(M)$, and choose an orientation-preserving 
embedding $\phi \from \D^n \to M$ with $B = \phi(\D^n)$. The derivative~$d\phi$
is homotopic to a trivialization $\wt{d\phi} \from \D^n \times \R^n \to TM \res
B$ (covering~$\phi$) in which the hyperplane field $\xi \res B$ is the kernel of
the one-form $\wt{d\phi}_* dx_n$. Any map $g \from \D^n \to \S^{n-1}$ that is 
constant equal to $(0, \dots, 0, 1)$ near $\partial \D^n$ can then be used to 
construct a new hyperplane field $\eta = (g \cdot \xi)_{\wt{d\phi}}$ on~$M$: set
$\eta = \xi$ out of~$B$ and, regarding $\S^{n-1}$ as the unit sphere of the dual 
space $(\R^n)^*$, define $\eta(p)$ for $p \in B$ to be the kernel of the linear 
form $\wt{d\phi}_* g(p)$.

This construction induces a group action of $\pi_n \S^{n-1}$ on $\pi_0 
\PF(M)$ and $\pi_0 \PF_c(M)$. First observe that, since $M$ is connected as well
as the implied spaces of embeddings and trivializations, the homotopy class of 
the hyperplane field $(g \cdot \xi)_{\wt{d\phi}}$ is unsensitive to the choice 
of $\phi, \wt{d\phi}$ and depends only on the homotopy classes of $\xi$
and~$g$.
Next, the constant map $g = (0, \dots, 0, 1)$ acts trivially. Finally, to check
the composition rule, denote by $\wt{d\phi}_\pm$ the restriction of $\wt{d\phi}$
to $\D^n_\pm \times \R^n$ where
$$\D^n_\pm = \{ (x-1, \dots, x_n) \in \D^n \mid 
\pm x_n \ge 0 \}.$$
If $g \from \D^n_- \to \S^{n-1}$ and $h \from \D^n_+ \to 
\S^{n-1}$ are maps equal to $(0, \dots, 0, 1)$ near the boundary and if $gh$ 
refers to the resulting map $\D^n \to \S^{n-1}$ ---~whose homotopy class is the 
product of the homotopy classes of $g$ and~$h$~---, then
$$ (gh \cdot \xi)_{\wt{d\phi}} =
   \bigl( g \cdot (h \cdot \xi)_{\wt{d\phi}_+} \bigr)_{\wt{d\phi}_-} \,. $$
This action and the obstruction cocycle determine the homotopy classification of
hyperplane fields on~$M$:

\begin{proposition} \label{p:homotopy}
Let $M$ be a connected, oriented manifold of dimension $n \ge 2$ with zero Euler
characteristic. The continuous cocycle 
$$ c \from \PF_c(M) \times \PF_c(M) \longrightarrow H_1 (M;\Z) $$
is surjective and two hyperplane fields $\xi, \eta$ satisfy $c(\xi,\eta) = 0$ if
and only if their homotopy classes are in the same orbit of $\pi_n \S^{n-1}$.
Furthermore, the stabilizer in $\pi_n \S^{n-1}$ of the homotopy class of any
hyperplane field~$\xi$ is trivial if $n \ge 4$ and, for $n=3$, is the image of
the homomorphism $H_2(M;\Z) \to \pi_1 \SO_2 = \pi_3 \S^2$ defined by the Euler 
class of~$\xi$.
\end{proposition}

Strictly speaking, as an obstruction class, the Euler class of~$\xi$ belongs to
$H^2 (M, \pi_1 \S^1)$, and so the implied homomorphism is rather given by the 
$\pi_1 \SO_2$--valued lift of the second Stiefel--Whitney class of~$\xi$. On the 
other hand, $\pi_1 \SO_2$ is identified with $\pi_3 \S^2$ via the Pontryagin 
isomorphism (see \fullref{x:pontryagin} below).

The proof of this proposition is actually more instructive than its statement.
We first recall a simple fact: for any finite dimensional vector space~$E$, the
tangent space to the Grassmann manifold $\G_k(E)$ at each point $\tau$ (a vector
subspace of dimension~$k$ in~$E$) can be canonically identified with $\tau^\ast
= \Hom (\tau, E/\tau)$. With this in mind, the key observation is the following:

\begin{lemma} \label{l:cobordism}
Let $M$ be a connected, oriented manifold of dimension $n \ge 2$ with zero Euler
characteristic. Given a hyperplane field $\xi \in \PF_c(M)$, there is a natural
one-to-one correspondence between $\pi_0 \PF_c(M)$ and the set $\Omega_1 (M;
\xi^\ast)$ of cobordism classes of $\xi^\ast$--framed curves in~$M$.
\end{lemma}

A \emph{$\xi^\ast$--framed curve} is a pair $(C, \gamma)$ consisting of a closed
one-dimensional submanifold~$C$ in~$M$ and a bundle equivalence (that is, a
bundle isomorphism over the identity) $\gamma \from \nu C \to \xi^\ast \res C$,
where $\nu C$ denotes the normal bundle of~$C$ and $\xi^\ast$ the bundle $\Hom 
(\xi, TM/\xi)$ ---~which can be identified with the dual bundle $\xi^*$ if $\xi$
is given as the kernel of a one-form, or with $\xi$ itself if $M$ is equipped 
with a metric. Two $\xi^\ast$--framed curves $(C, \gamma)$ and $(C', \gamma')$ 
are \emph{cobordant} if there exists a compact surface~$S$ in $[0,1] \times M$, 
with $\partial S = (\{0\} \times C) \cup (\{1\} \times C')$, such that $\gamma 
\cup \gamma' \from \nu S \res{\partial S} \to \xi^\ast \res{\partial S}$ extends
to a bundle equivalence $\nu S \to \xi^\ast \res S$, where indeed $\xi^\ast$ 
stands here for its pullback over $[0,1] \times M$ (we will consistently use the
same notation for a bundle over~$M$ and its pullback over $[0,1] \times M$).

\begin{proof}
From hyperplane fields to $\xi^\ast$--framed curves, the correspondence goes as
follows. For a generic hyperplane field $\eta \in \PF_c(M)$, the set $C(\xi,\eta
)$ of points where $\xi$ coincides with~$-\eta$ is a $\xi^\ast$--framed curve. In
fact, if we regard $\xi$ and~$\eta$ as sections of the sphere cotangent bundle 
$ST^*M$ and denote their images by $X$ and~$Y$, respectively, $C = C(\xi,\eta)$
is the projection to~$M$ of the intersection $X \cap (-Y)$. If $X$ and~$(-Y)$ 
are transverse to each other, $C$ is a closed curve. Furthermore, over any point
$p$ in~$M$, the tangent space of~$X$ determines a projection from the tangent 
space of $ST^*M$ at $\xi(p)$ to $\xi^\ast(p)$, the tangent space of the fiber 
$ST^*_pM$. Along the curve~$C$, since transversality holds, the composition of 
the differential of~$-\eta$ with this projection provides the required bundle 
equivalence $\gamma \from \nu C \to \xi^\ast \res C$. Clearly, the cobordism 
class of the $\xi^\ast$--framed curve $(C, \gamma)$ only depends on the homotopy
class of~$\eta$.

The correspondence in the other direction is a version of the Thom--Pontryagin
construction. Given a $\xi^\ast$--framed curve $(C, \gamma)$, the considerations
above show how $\gamma$ defines a germ~$\eta_0$ of hyperplane field near~$C$
(actually only the $1$--jet of it) that coincides transversally with $-\xi$ along
$C$. Then pick one-forms $\alpha$ and~$\beta_0$ defining $\xi$ and $\eta_0$, 
respectively, choose a function $\rho \from M \to [0,1]$ equal to~$1$ in the 
$\eps$--neighborhood of~$C$ and to~$0$ out of the $2\eps$--neighborhood of~$C$, 
and consider the one-form $\beta = (1-\rho) \alpha + \rho \beta-0$. For $\eps$ 
sufficiently small, $\beta$~is defined everywhere and non-singular, and its 
kernel~$\eta$ is a hyperplane field whose associated $\xi^\ast$--framed curve is
$(C,\gamma)$. The parametric version of this construction associates a homotopy
of hyperplane fields to any cobordism between $\xi^\ast$--framed curves, so the 
correspondence between homotopy classes and cobordism classes is well-defined 
and one-to-one.
\end{proof}

\begin{example}[The Pontryagin isomorphism] \label{x:pontryagin}
Let $\xi$ be a hyperplane field on $B = \D^n$. Choose a positive trivialization
of $TB$ in which $\xi$ is spanned by the first $(n-1)$ basis vectors, and note 
that $\xi$ and $\xi^\ast$ are then also trivialized. Thus, $\xi^\ast$--framed 
curves in~$B$ are just usual framed curves while elements of $\PF_c (B)$ are 
identified with maps $B = \D^n \to \S^{n-1}$ that are constant equal to $(0, 
\dots, 0, 1)$ near $\partial B$. So we get a natural one-to-one correspondence 
$$ \pi_n \S^{n-1} = \pi_0 \PF_c (B) \longrightarrow
   \Omega_1 (B; \xi^\ast) = \Omega_1 (B; \ul\R^{n-1}) $$
which is indeed independent of $\xi$ and of the trivialization since it takes
the homotopy class of any generic $g \from \D^n \to \S^{n-1}$ to the cobordism
class of the framed curve $(C_g, \gamma_g)$ defined as follows:
\begin{itemize}
\item
$C_g$ is the fiber $g^{-1}(q)$ where $q = (0, \dots, 0, -1) \in \S^{n-1}$ is 
supposed to be a regular value;
\item
$\gamma_g \from \nu C_g \to \xi^\ast \res{ C_g } = C_g \times T_{-q}\S^{n-1}$ is
given by $\gamma_g (p,w) = \bigl( p, -dg(p)\, w \bigr)$.
\end{itemize}
Furthermore, this correspondence induces a group structure on $\Omega_1 (B ; \ul
\R^{n-1})$: any two framed curves can be individually isotoped into disjoint 
balls and the sum of their cobordism classes is then the cobordism class of the
union.

Now let $C \subset \Int B$ be a circle in the $x_1x_2$--plane and let $\gamma^1 
\from \nu C \to C \times \R^{n-1}$ be its standard normal framing ---~the normal 
vector in the plane followed by the canonical basis of the remaining $\R^{n-2}$.
Any loop $u \from \S^1 = C \to \SO_{n-1}$, considered as an automorphism of $C 
\times \R^{n-1}$, can be composed with $\gamma^1$ to give a framing $\gamma^u = 
u \cdot \gamma^1$. The resulting map 
$$ \pi_1 \SO_{n-1} \longrightarrow 
   \Omega_1 (B; \ul\R^{n-1}) = \pi_n \S^{n-1} $$
is a group homomorphism (for, if copies of $(C, \gamma^u)$ and $(C, \gamma^v)$
are placed in disjoint balls, their union is cobordant to $(C, \gamma^{uv})$) 
and is surjective (because any cooriented closed curve in~$B$ is cobordant to 
$C$). By a theorem of L~Pontryagin, this map is indeed an isomorphism, and the 
arguments in the proof below will actually show that it is injective.
\end{example}

\begin{proof}[Proof of \fullref{p:homotopy}]
Note first that the orientations of $\xi$ and~$M$ induce an orientation of~$\xi
^\ast$. Therefore, any $\xi^\ast$--framed curve $(C, \gamma)$ is (co)\,oreiented
by $\gamma \from \nu C \to \xi^\ast \res C$, and so it determines a $1$--cycle 
whose homology class depends only on the cobordism class of $(C, \gamma)$. This
gives a map $\Omega_1 (M;\xi^\ast) \to H_1 (M;\Z)$ whose composition with the
bijection $\pi_0 \PF_c(M) \to \Omega_1 (M;\xi^\ast)$ sends the homotopy class of
a hyperplane field~$\eta$ to the obstruction class $c(\xi,\eta)$. The cocycle 
$c$ is then surjective for, in dimension $n \ge 2$, any homology class in $H_1 
(M;\Z)$ can be represented by an embedded, oriented, closed curve, and any such
curve~$C$ admits a $\xi^\ast$--framing since both $\nu C$ and $\xi^\ast \res C$ 
are trivial bundles.

Consider two hyperplane fields $\xi, \eta \in \PF_c(M)$. By definition, their 
homotopy classes lie in the same $\pi_n \S^{n-1}$--orbit if and only if $\eta$ is
homotopic to a hyperplane field that coincides with~$\xi$ out of a ball, and 
this property clearly implies $c(\xi,\eta) = 0$. Conversely, if $c(\xi,\eta)$ is
zero and if $\xi$ and~$\eta$ are generic, the curve $C = C(\xi,\eta)$ consisting
of points where $\xi = -\eta$ is nullhomologous. Hence, $C$~is cobordant to $C'
= \partial D$ where $D$ is an embedded disk disjoint from~$C$. Let $S \subset 
[0,1] \times M$ be a connected, oriented, compact surface with $\partial S = 
(\{0\} \times C) \cup (\{1\} \times C')$. Since $S$ retracts onto the union of 
$\{0\} \times C$ and a graph, the $\xi^\ast$--framing of~$C$ extends over~$S$. 
Therefore, by \fullref{l:cobordism}, $\eta$~is homotopic to a hyperplane field
$\eta'$ such that $C(\xi,\eta') = C'$, and using an affine homotopy away from
$C'$, we can indeed make~$\eta'$ equal to~$\xi$ out of a small neighborhood of 
$C' = \partial D$. Thus, the homotopy classes of $\xi$ and~$\eta$ are in the 
same orbit under $\pi_n \S^{n-1}$.

It remains to compute the stabilizer of the homotopy class of~$\xi$. To do this,
take a ball $B \subset M$ parameterized by~$\D^n$. As explained in
\fullref{x:pontryagin}, the cobordism classes corresponding via \fullref{l:cobordism}
to the $\pi_n \S^{n-1}$--orbit of the homotopy class of~$\xi$ are represented by
the $\xi^\ast$--framed curves $(C, \gamma^u)$, where $C$ is a fixed circle in~$B$
and $u$ a loop in $\SO_{n-1}$. Our task is to determine the loop classes $[u] 
\in \pi_1 \SO_{n-1}$ for which the $\xi^\ast$--framing $\gamma^u$ extends over 
some connected, oriented, compact surface in $[0,1] \times M$ bounded by $\{0\} 
\times C$. 

Let $S$ be such a surface and choose a trivialization of the bundle $\nu S$ 
---~and thereby also a trivialization of $\nu C$. This choice identifies bundle
equivalences $\nu S \to \xi^\ast \res S$ with trivializations of $\xi^\ast \res
S$. But any two trivializations of this bundle are  homotopic over the boundary
because they differ by a map $S \to \SO_{n-1}$, and the restriction of such a 
map to $\partial S = \{0\} \times C$ (which is connected) is nullhomotopic. This
shows that, up to homotopy, there exists a unique $\xi^\ast$--framing of~$C$ that
extends over the given surface~$S$.

Consider now a disk $D \subset [-1,0] \times M$ bounded by $\{0\} \times C$ and
whose projection to~$M$ is the affine disk spanned by $C$ in~$B$. The (unique) 
$\xi^\ast$--framing of $C$ that extends to~$D$ is the standard framing~$\gamma^1$
(see \fullref{x:pontryagin}). Denote by $\ol S \subset \R \times M$ the 
closed, connected, oriented surface obtained by smoothing $S \cup (-D)$ in the 
obvious way. Since $T(\R \times M) = TM \oplus \ul\R = \xi \oplus \ul\R^2$ while
$T(\R \times M) \res{\ol S} = \nu\ol S \oplus T\ol S$ ---~and $T\ol S \oplus \ul
\R = \ul\R^3$, the bundles $\xi^\ast \res{ \ol S }$ and $\nu\ol S$ are stably 
equivalent. If $n \ge 4$ then $\xi^\ast \res{\ol S}$ and $\nu\ol S$ are indeed 
equivalent (because a vector bundle of rank at least four over a surface has a
connected space of non-vanishing sections), and so the $\xi^\ast$--framing of~$C$
that extends over~$D$ also extends over~$S$. Therefore, $\gamma^1$ is the unique
$\xi^\ast$--framing of~$C$ that is nullcobordant, and so $\pi_n\S^{n-1} = \Z/2\Z$
acts freely. If $n=3$, however, the bundle $\nu\ol S$ is trivial but $\xi^\ast 
\res{\ol S}$ is not in general: given trivializations of $\xi^\ast \res D$ and 
$\xi^\ast \res S$, the induced trivializations of $\xi^\ast \res C$ differ by a
map $C \to \SO_2$ whose degree is the Euler class $e(\xi)$ of~$\xi$ evaluated on
$[\ol S] \in H_2(M;\Z)$. Thus, the $\xi^\ast$--framing $\gamma^u$ of~$C$ extends
over~$S$ if and only if $\langle e(\xi), [\ol S] \rangle = [u] \in \pi_1 \SO_2$.
This completes the proof.
\end{proof}

Let's now return to our three-dimensional framework. It follows from \fullref{p:homotopy} that the stabilizer of the homotopy class of a plane field
$\xi$ under the action of $\pi_3 \S^2 = \Z$ is $\abs \xi \Z$, where $\abs \xi$
denotes the divisibility of the (torsion-free part of the) Euler class of~$\xi$.
For any plane field~$\eta$ homologous to~$\xi$, we will call \emph{relative
framing} of $\xi$ and~$\eta$ the element $d(\xi,\eta) \in \{ 0, \dots, \abs \xi
- 1 \}$ which takes the homotopy class of $\xi$ to that of~$\eta$. As an 
illustration, we can recover a calculation of W~Neumann and L~Rudolph:

\begin{lemma}[Neumann--Rudolph \cite{NR}] \label{l:framing}
Let $\xi$ denote the standard contact structure on~$\S^3$ ---~that is, the plane
field orthogonal to the Hopf fibers~--- and $\xi^-$ the plane field associated
with the negative Hopf link~$H^-$. Then the relative framing $d(\xi,\xi^-)$ is
equal to~$1$.
\end{lemma}

\begin{proof}
Since the Hopf flow preserves the open book given by $H^-$ (whose mapping to the
circle is the argument of the Hopf fibration: see \fullref{x:plumbing}), we
can arrange that it preserves $\xi^-$ too (just construct $\xi^-$ as indicated
in the introduction). On the other hand, $\xi$~is the plane field orthogonal to
the Hopf fibers, and is also invariant under the Hopf flow. Therefore, the
$\xi^\ast$--framed curve $(C, \gamma)$ determined by $\xi^-$ is also invariant.
Clearly, the curve~$C$ (the set of points where $\xi^-$ coincides with $-\xi$)
is the component $\{ r_2 = 0 \}$ of~$H^-$, which is transverse to~$\xi$. Then,
if we identify $\xi^\ast$ with~$\xi$ using the metric of~$\S^3$, the bundle
equivalence $\gamma \from \nu C \to \xi \res C = \nu C$ is (homotopic to) the
identity.

Now let $D$ be a disk bounded by~$C$ in $[0,1] \times \S^3$. The trivialization
of $\nu C$ that extends to $\nu D$ differs from the trivialization of $\xi \res
C$ that extends to $\xi \res D$ by one twist: indeed, $\xi$~admits a global
non-vanishing section and the linking number of~$C$ and its push-off along this
section is equal to~$-1$. This proves that the relative framing $d(\xi,\xi^-)$
is equal to~$1$.
\end{proof}

\section{Contact structures}
\label{sec:c}

We briefly present here some notions and results of three-dimensional contact
geometry that we will invoke to prove \fullref{t:equivalence}. A \emph
{contact form} on an oriented three-manifold~$M$ is a one-form~$\alpha$ whose
exterior product with $d\alpha$ is everywhere positive ---~with respect to the
orientation of $M$. A \emph{contact structure} on $M$ is a (co)\,oriented plane
field~$\xi$ which is the kernel of some contact form, and a \emph{contact
manifold} is a manifold equipped with a contact structure. A fundamental
property of contact structures, established by J~Gray, is that they are
$\class1$--stable: if $\xi_s$, $s \in [0,1]$, is a path of contact structures
on a closed manifold~$M$, then there exists an isotopy $\phi_s$ of~$M$ such that
$\phi_0 = \id$ and $\phi_{s*} \xi_0 = \xi_s$ for all $s \in [0,1]$. Thus, two
contact structures on a closed manifold~$M$ are isotopic if and only if they are
in the same homotopy class of contact structures (that is, the same
connected component of the space of contact structures).

The possibility of constructing contact structures on three-manifolds from
open books was discovered by Thurston and Win\-kelnkemper~\cite{TW}.
However, the systematic study of the relations between these two geometric
objects is much more recent (see the article \cite{Gi} by the first
author) and is based on the following:

\begin{definition} \label{d:compatibility}
Let $M$ be a closed oriented three-manifold. We say that a contact structure
$\xi$ on~$M$ \emph{is carried} by an open book $(K,\theta)$ if it is the kernel
of a one-form~$\alpha$ satisfying the following conditions:
\begin{itemize}
\item
$\alpha$ induces a positive non-singular form on~$K$;
\item
$d\alpha$ induces a positive area form on each fiber of~$\theta$.
\end{itemize}
Any such one-form $\alpha$ is said to be adapted to $(K,\theta)$.
\end{definition}

With this terminology, the main result of Thurston and
Winkelnkemper \cite{TW} is that the set of contact
structures carried by a given open book is non-empty. It is easy to check that
this set is also open and contractible in the space of all contact structures on
$M$ (see the proof of \fullref{l:overtwisted}). In particular, according to
Gray's stability theorem, all the contact structures it contains belong to the
same isotopy class. For this reason, we often speak of ``the contact structure
associated with the open book'', this contact structure being defined only up to
isotopy.

\begin{remark} \label{r:compatibility}
The contact structure associated with an open book $(K,\theta)$ ---~in a closed
oriented three-manifold~$M$~--- belongs to the homotopy class of plane fields
associated with $(K,\theta)$. In fact, if $\alpha_0$ is a one-form defining the
plane field associated with $(K,\theta)$ as in the introduction (with $N$ small
enough) and if $\alpha_1$ is a contact form adapted to $(K,\theta)$, then all
forms $(1-t) \alpha_0 + t \alpha_1$, $t \in [0,1]$, are non-singular, and so
their kernels yield the desired homotopy of plane fields.
\end{remark}

\begin{example} \label{x:compatibility}
The standard contact structure $\xi$ on $\S^3 \subset \C^2$ is defined by the
one-form $\alpha = r_1^2 d\theta_1 + r_2^2 d\theta_2$. This contact form is
adapted to the trivial open book $(U, \theta_1)$ for it induces on each fiber
of~$\theta_1$ the one-form $r_2^2 d\theta_2$. It is also adapted to the open
book $(H^+, \theta_1 + \theta_2)$, for it induces on each fiber of $\theta_1 +
\theta_2$ the one-form $(2r_1^2-1) \, d\theta_1 = (2r_2^2-1) \, d\theta_2$.
\end{example}

The above example shows that a given contact structure may be carried by several
open books. Indeed, according to a previous article by the first author
\cite{Gi}, any contact structure is carried by
many open books but we have the following stable equivalence theorem:

\begin{theorem}[Giroux \cite{Gi}] \label{t:pequivalence}
On a closed oriented three-manifold, two open books carrying the same contact
structure admit isotopic positive stabilizations. 
\end{theorem}

The last ingredient of contact geometry we will need is the classification of
overtwisted contact structures, due to Y~Eliashberg. A contact structure~$\xi$
on a three-manifold~$M$ is \emph{overtwisted} if there exists a simple closed
curve $L \subset M$ with the following properties:
\begin{itemize}
\item
$L$ is Legendrian, that is, is tangent to~$\xi$ at each point;
\item
$L$ is unknotted, that is, bounds a disk;
\item
the Thurston--Bennequin number of~$L$ ---~that is, the linking number of~$L$
and its push-off along the normal vector to~$\xi$~--- is non-negative.
\end{itemize}
Overtwisted contact structures have an extremely simple classification:

\begin{theorem}[Eliashberg \cite{El}] \label{t:overtwisted}
On a closed oriented three-manifold, two overtwisted contact structures are
isotopic if and only if they are in the same homotopy class of plane fields.
\end{theorem}

\section{Proof of the stable equivalence theorem}
\label{sec:d}

Let $(K,\theta)$ be an open book in a closed oriented three-manifold~$M$, with
associated plane field~$\xi$, and let $(K^+, \theta^+)$ and $(K^-, \theta^-)$
denote open books obtained from $(K,\theta)$ by $H^+$--plumbing and
$H^-$--plumbing, respectively.

\begin{lemma} \label{l:homotopy}
The plane field $\xi^+$ associated with $(K^+, \theta^+)$ is homotopic to~$\xi$
while the plane field~$\xi^-$ associated with $(K^-, \theta^-)$ is homologous to
$\xi$ with relative framing $d(\xi,\xi^-)$ equal to~$1$.
\end{lemma}

\begin{proof}
Consider the open books in $\S^3$ given by $H^+$ and by the unknot~$U$. After
isotoping one of them, we may assume that they have a common plumbing ball~$B$
(see \fullref{sec:a}) in which they coincide as well as their associated plane fields.
Now these plane fields are homotopic to each other (according to \fullref
{x:compatibility} and \fullref{r:compatibility}, both are homotopic to the
standard contact structure) and, since $\S^2$ is simply connected, they are also
homotopic relative to~$B$. Therefore, $\xi^+$ is homotopic to the plane field
associated with the open book obtained by plumbing $(K,\theta)$ with~$U$. But
this open book is isotopic to $(K,\theta)$, so $\xi^+$ is homotopic to~$\xi$. On
the other hand, $\xi^-$ and~$\xi$ are homologous since they coincide out of a
ball. Next, arguing as above, we see that the relative framing $d(\xi,\xi^-)$ is
equal to the relative framing in~$\S^3$ of the standard contact structure and
the plane field associated with the negative Hopf link~$H^-$. Then the result
follows from \fullref{l:framing}.
\end{proof}

The next lemma is essentially due to Torisu:

\begin{lemma}[Torisu \cite{To}] \label{l:overtwisted}
The contact structure $\xi^-$ associated with $(K^-,\theta^-)$ is overtwisted.
\end{lemma}

\begin{proof}
Let's say that a one-form $\beta$ on a compact oriented surface is admissible if
it induces a positive non-singular form on the boundary and if its differential
$d\beta$ is a positive area form in the interior. A contact form adapted to an
open book clearly induces an admissible one-form on each page. The existence and
uniqueness ---~up to isotopy~--- of contact structures carried by an open book
is mostly due to the contractibility of the space of admissible forms on a given
surface. This contractibility also allows to construct an adapted contact form 
inducing a prescribed admissible form on a given page (Thurston--Winkelnkemper
\cite{TW}).

Now let $F^- = F \cup A^-$ be a page of $(K^-,\theta^-)$, where $F$ is a page of
$(K, \theta)$ and~$A^-$ a negative Hopf band. Since the core curve $L$ of~$A^-$
is homologically non-zero in~$F^-$, there exists an admissible form~$\beta^-$
that vanishes at each point of~$L$. Then consider a contact form~$\alpha^-$ on
$M$ which is adapted to $(K^-,\theta^-)$ and induces $\beta^-$ on~$F^-$. For the
contact structure $\xi^-$ defined by~$\alpha^-$, the curve~$L$ is Legendrian and
unknotted. Moreover, the normal vector to $\xi^-$ along~$L$ is the normal vector
to~$A^-$, so the Thurston-Bennequin number of~$L$ is equal to~$1$. Thus, the
contact structure~$\xi^-$ is overtwisted.
\end{proof}

\begin{proof}[Proof of \fullref{t:equivalence}]
If two open books admit isotopic stabilizations, their associated plane fields
are homologous since the homology class does not change under $H^\pm$--plumbing.
Suppose now that $(K,\theta)$ and $(K',\theta')$ are two open books in~$M$ whose
associated contact structures $\xi$ and $\xi'$ are homologous as plane fields.
Let $d$ denote the relative framing $d(\xi,\xi')$ and consider an open book
$(K'',\theta'')$ obtained from $(K,\theta)$ by $d$~successive $H^-$--plumbings.
According to \fullref{l:homotopy} and \fullref{l:overtwisted}, the contact
structure~$\xi''$ associated with $(K'',\theta'')$ is in the same homotopy class
of plane fields as $\xi'$ and is overtwisted provided $d \ge 1$. Applying one
more $H^-$--plumbing to both $(K',\theta')$ and $(K'',\theta'')$ if necessary, we
may assume that $\xi'$ and~$\xi''$ are both overtwisted. Then it follows from
Eliashberg's \fullref{t:overtwisted} that $\xi'$ and $\xi''$ are actually
isotopic. Hence \fullref{t:pequivalence} implies that $(K',\theta')$ and
$(K'',\theta'')$ admit isotopic positive stabilizations, so $(K,\theta)$ and
$(K',\theta')$ admit isotopic stabilizations.
\end{proof}

\begin{remark}
The above proof shows that, while there is no control on the number of necessary
$H^+$--plumbings, the number of necessary $H^-$--plumbings can be bounded \emph{a
priori} in terms of the relative framing of $\xi$ and $\xi'$, namely,  by 
$$2 + \min \{ d(\xi,\xi'), d(\xi',\xi) \}.$$
\end{remark}

Corollaries \ref{c:links} and \ref{c:harer} follow readily from
\fullref{t:equivalence}. \fullref{c:knots} follows similarly from
the following
refined version of \fullref{t:equivalence}:

\begin{theorem}
In a closed oriented three-manifold, two open books with connected bindings and
homologous associated plane fields admit isotopic stabilizations which can be
obtained by finitely many successive $T^+$--plumbings and $E$--plumbings.
\end{theorem}

\begin{proof}
Suppose that $(K,\theta)$ and $(K',\theta')$ are two open books in~$M$ whose
binding are connected and whose associated contact structures $\xi$ and $\xi'$ 
are homologous as plane fields. Since $E$--plumbing is a composition of an $H^+$--
and an $H^-$--plumbing, its effect on the homotopy class and the isotopy type of
the associated contact structure is the same as the effect of $H^-$--plumbing. 
Therefore, after performing a number of $E$--plumbings on our open books (see the
proof of \fullref{t:equivalence}, we may assume that $\xi$ and $\xi'$ are
overtwisted and homotopic as plane fields, and hence isotopic by Eliashberg's
\fullref{t:overtwisted}. Then we conclude with the following refined version
of \fullref{t:pequivalence}: on a closed oriented three-manifold, two open 
books carrying the same contact structure and having connected bindings admit 
isotopic positive stabilizations which can be obtained by $T^+$--plumbings. 
\end{proof}

\begin{proof}[Proof of \fullref{c:grothendieck}]
\fullref{c:harer} shows that the Grothendieck group $\Gamma$ of fibered
links in~$\S^3$ is generated by the Hopf links $H^+$ and~$H^-$. To complete the
proof, we proceed as Neumann--Rudolph in~\cite{NR}. To each fibered
link $K$ in~$\S^3$ we assign two integers, $\mu(K)$ and $\lambda(K)$:
\begin{itemize}
\item
$\mu(K)$ is the Milnor number of~$K$, that is, the first Betti number  of a
fiber Seifert surface (a page of the corresponding open book);
\item
$\lambda(K)$ is the ``enhanced Milnor number'', that is, relative framing of
the plane field associated with the unknot (the standard contact structure for
instance) and the plane field associated with~$K$.
\end{itemize}
The additivity of $\mu$ and~$\lambda$ under plumbing (which follows from our
discussion in Sections \ref{sec:a} and \ref{sec:b} ---~see also the proof of \fullref{l:homotopy})
implies that the pair $(\mu,\lambda)$ induces a homomorphism from $\Gamma$ to
$\Z^2$. By \fullref{l:framing} this homomorphism maps the generators $H^+$ and
$H^-$ to $(1,0)$ and $(1,1)$, respectively, and so it is an isomorphism.

The calculation of the Grothendieck group of fibered knots is analogous.
\end{proof}

\bibliographystyle{gtart}
\bibliography{link}
\end{document}